\documentclass[12pt]{article}
\title{A new approach to flatness in the category of quasi-coherent sheaves}
\date{}
\author{}
\usepackage{graphicx,fancyhdr,fancybox,rotating,xcolor,appendix}
\usepackage{hhline,multicol}
\usepackage{amsthm,amsbsy,geometry,times}
\usepackage{amsmath,tikz-cd}
\usepackage{lipsum}
\usepackage{euscript}
\usepackage{hyperref}
\usepackage{hyperref}
\usepackage{hyperref}
\usepackage{pdfpages} 
\usepackage{textcomp}
\usepackage{stmaryrd}
\usepackage[all]{xy}
\usepackage{amsfonts}
\usepackage{amsmath}
\usepackage{amssymb}
\usepackage{latexsym}
\usepackage{mathrsfs}

\newtheorem{thm}{Theorem}[section]
 \newtheorem{cor}[thm]{Corollary}
 \newtheorem{lem}[thm]{Lemma}
 \newtheorem{prop}[thm]{Proposition}
 \newtheorem{Def}[thm]{Definition}
\newtheorem{rem}[thm]{Remark}

\def\hom{{\rm {\mathcal{H}om}_{\O}}}

\def\hm{{\rm {\mathcal{H}om}^{qc}_{\mathcal{O}_{X}}}}

\def\Ker{{\rm Ker}}

\newcommand{\X}{\O$-$Mod}
\newcommand{\XU}{\O(U)$-$Mod}
\newcommand{\U}{\mathcal{O}_U$-$Mod}
\newcommand{\Ui}{\mathcal{O}_{U_i}$-$Mod}
\newcommand{\Xx}{\mathcal{O}_{X,x}$-$Mod}
\newcommand{\Q}{\mathfrak{Qcoh}}

\usepackage[all]{xy}

\def\spec{{\rm Spec}}

\def\inj{{\rm \mathcal{I}nj}}

\def\Ext{{\rm Ext}}

\def\Hom{{\rm Hom}}

\def\Ker{{\rm Ker}}
\def\Flat{{\rm Flat}}
\def\Abs{{\rm Abs}}

\def\D{{\rm D}}

\def\O{{\mathcal{O}_{X}}}

\def\ApQ{{\mathcal{A}p}_{\Q(X)}}

\def\Apf{{{\mathcal{A}p}_{\Q(D_{A}(f_i))}}}

\def\Apxx{{{\mathcal{A}p}_{\mathcal{O}_{X,x}-Mod}}}

\def\Apq{{{\mathcal{A}p}_{\Q(X)}}}
\def\LAp{{{\mathcal{LA}p}_{\Q(X)}}}
\def\Apou{{{\mathcal{A}p}_{\mathcal{O}_X(U)-Mod}}}
\def\Apoui{{{\mathcal{A}p}_{\mathcal{O}_X(U_i)-Mod}}}

\begin{document}

\thispagestyle{empty}

\maketitle \vspace*{-1.5cm}
\begin{center}{\large\bf   Soumia Mamdouhi}
\bigskip

\small{ Moroccan School of Engineering Sciences (EMSI), Rabat, Morocco\\

  soumia.mamdouhi86@gmail.com
}

\end{center}

\noindent{\large\bf Abstract.} 
Recently, many authors have studied  certain properties of modules
such as projectivity, injectivity and flatness from an alternative point of view. This way, Durğun has introduced absolutely pure domains of modules as a mean to gauge the extent of the flatness of any given module (so not just to determine whether or not the module is flat). In this paper, we develop a new treatment of the goemetric flatness  in the category of quasi-coherent sheaves which shed more light on some of its various important aspects. Namely, in terms of this new approach, some classical schemes are characterized, some new results are established and also known ones  are generalized.
\bigskip


\bigskip\bigskip
%

\bigskip


%


\section{Introduction}
Throughout this paper, all rings are commutative with identity,  $X$ will denote a scheme. We denote by $\X$ and $\Q(X)$ the category of  $\O$-modules   and  the category of  quasi-coherent $\O$-modules, respectively. For an  open subset $U$  of $X$, we use  $\mathcal{O}_{U}$ to denote  the restriction of $\O$ to $U$, and we denote by $\U$  the category of $\mathcal{O}_U$-modules. As usual, for a point $x$ in $X$, we denote by $\Xx$ the category of $\mathcal{O}_{X,x}$-modules.
\medskip

\par The study of flatness and purity occupies a central place in  the category of 
quasi-coherent sheaves $\Q(X)$ ( see  \cite{PIS} and \cite{PMONOIDAL}). In particular, Enochs, 
Estrada and Odaba\c{s}\i\: proved in \cite{PIS} that pure injective and 
absolutely pure objects in $\Q(X)$ behave well with respect to the 
geometric structure of $X$. More recently, Aghasi and Nemati studied in 
\cite{absolute} absolute purity in $\Q(X)$ and established several 
characterizations of absolutely pure quasi-coherent sheaves.

\medskip

\par  On the other hand, many  studies are done every year on flat,  projective and injective  modules. Many of
them involve concepts derived from relative flatness,  projectivity and  injectivity. Rather
than saying whether a module has a certain property or not, each module is assigned a
relative domain that, somehow, measures to which extent it has this particular property.

\medskip

\par In \cite{Dur}, an alternative perspective on the flatness of a module was introduced. Recall that, for two modules $M$ and  $N$, $N$ is said to be \textit{absolutely $M$-pure} if the canonical 
map $M\otimes N\rightarrow M\otimes K$ is a monomorphism for every 
extension $K$ of $N$. The \textit{absolutely pure domain} of $M$ is then 
defined as the collection of all absolutely $M$-pure modules. Note that 
$M$ is flat if and only if its absolutely pure domain is the entire module 
category, while at the other extreme, $M$ is called \textit{$f$-indigent} 
if its absolutely pure domain is minimal, consisting exactly of the 
absolutely pure modules. In \cite{Dur}, properties of absolutely pure 
domains and $f$-indigent modules are studied in detail, the existence of 
$f$-indigent modules is proved for arbitrary rings, and several classes of 
rings are characterized via these notions.

\medskip

\par In this  paper, we extend the theory developed in \cite{Dur} to
the category $\Q(X)$. This study provides new results that contribute to the study of geometric flatness and purity in $\Q(X)$.
 
\medskip

 The paper is organized as follows:\\
\medskip
In Section 2, we recall the essential concepts necessary for this article.\\
\medskip
In Section 3, we introduce absolutely pure domains in $\Q(X)$: for 
$\mathcal{M}\in\Q(X)$, we write $\ApQ(\mathcal{M})$ for the class of 
objects that are absolutely $\mathcal{M}$-pure, the direct analogue of 
absolutely pure domain  in \cite{Dur}. Our first main result describes the absolutely pure domains of some objects in $\Q(X)$
 in terms of Ext-orthogonal classes. Recall that for an object 
$\mathcal{M}\in\Q(X)$, its left Ext-orthogonal class is $^{\perp}{\mathcal{M}}=\{\mathcal{N}\in\Q(X)\mid \Ext^1(\mathcal{N},\mathcal{M})=0\}$.

\medskip
\noindent\textbf{Theorem A} (= Theorem \ref{M-F-M})\textbf{.} 
\textit{Let $0\to \mathcal{M} \to \mathcal{F}\to \mathcal{M}'\to 0$ be 
a short exact sequence in $\Q(X)$ with $\mathcal{F}$ flat. Then 
$^{\perp}{\mathcal{M}'^{+}}\subseteq \ApQ(\mathcal{M})$. If moreover 
$\inj(X)\subseteq {^{\perp}\mathcal{M}'^{+}}$, then 
${^{\perp}\mathcal{M}'^{+}}= \ApQ(\mathcal{M})$.}

\medskip
In general, the notion of absolutely pure domain is not local in $\Q(X)$, the condition $\mathcal{N}\in\ApQ(\mathcal{M})$ cannot be
checked on affine open subsets of $X$ (see Remark \ref{AA}). We
therefore introduce the \emph{locally absolutely pure domain}
$\LAp(\mathcal{M})$: an object $\mathcal{N}$ belongs to
$\LAp(\mathcal{M})$ if $\mathcal{N}(U)$ is absolutely
$\mathcal{M}(U)$-pure as $\mathcal{O}_X(U)$-modules, for every affine
open subset $U\subseteq X$. We show that
\[
\LAp(\mathcal{M})\subseteq\ApQ(\mathcal{M})
\]
for every $\mathcal{M}\in\Q(X)$, and that equality holds whenever $X$ is
locally coherent, semi-separated, and quasi-compact.

\medskip
\noindent\textbf{Theorem B} (= Theorem \ref{help coh})\textbf{.} 
\textit{Let $X$ be a locally coherent scheme. Then 
$\LAp(\mathcal{M}) \subseteq \ApQ(\mathcal{M})$ for every 
$\mathcal{M}\in\Q(X)$. If moreover $X$ is semi-separated and 
quasi-compact, then $\LAp(\mathcal{M}) = \ApQ(\mathcal{M})$.}

\medskip
Theorem B leads to the following characterization of locally coherent
schemes.

\medskip
\noindent\textbf{Theorem C} (= Theorem \ref{coh})\textbf{.} 
\textit{Let $X$ be a semi-separated quasi-compact scheme. Then $X$ is 
locally coherent if and only if $\ApQ(\mathcal{M}) \subseteq 
\LAp(\mathcal{M})$ for every $\mathcal{M}\in\Q(X)$.}

\medskip
In Section 4, we introduce the notion of $f$-indigent objects in $\Q(X)$: for $\mathcal{M}\in\Q(X)$, $\mathcal{M}$ is called $f$-indigent if its absolutely pure domain consists only of absolutely pure objects in $\Q(X)$, the direct analogue notion introduced in \cite{Dur}. In this section, we
first prove the existence of $f$-indigent objects when $X$ is
quasi-compact and quasi-separated, extending
\cite[Proposition 3.1]{Dur} to the category $\Q(X)$. We then introduce
the notion of \emph{locally $f$-indigent} objects: $\mathcal{M}\in\Q(X)$
is locally $f$-indigent if $\mathcal{M}(U)$ is $f$-indigent as an
$\mathcal{O}_X(U)$-module for every affine open subset $U\subseteq X$.
We establish their basic properties and show, in particular, that every
locally $f$-indigent object is $f$-indigent whenever $X$ is
semi-separated and quasi-compact.

As an application, we obtain a homological characterization of
Noetherian schemes in terms of pure injective indigent objects. More
precisely, we prove the following result, which extends
\cite{Dur} to the geometric setting.

\medskip
\noindent\textbf{Theorem D} (= Theorem \ref{noetherian})\textbf{.} 
\textit{Let $X$ be a quasi-compact locally coherent scheme. Then $X$ is 
Noetherian if and only if there exists a pure injective indigent 
quasi-coherent $\mathcal{O}_U$-module for every affine open $U\subseteq X$.}

\section{Preliminaries}

\par In this section, we recall the basic notions of purity, flatness 
and the internal hom functor in $\Q(X)$ that will be used throughout 
the paper.

\par Recall that $\Q(X)$ is a monoidal category with the tensor product 
\emph{induced} from $\mathcal{O}_X\text{-Mod}$: for $\mathcal{F}, \mathcal{G} \in \Q(X)$, 
$\mathcal{F}\otimes_{\mathcal{O}_X} \mathcal{G}$ is the sheafification of the presheaf 
$U\mapsto \mathcal{F}(U)\otimes_{\mathcal{O}_X(U)}\mathcal{G}(U)$, and this sheaf is again 
quasi-coherent, so $\Q(X)$ is closed under $\otimes_{\mathcal{O}_X}$. 
In contrast, the internal hom functor $\mathcal{H}om_{\mathcal{O}_X}(-,-)$ in $\X$
does \emph{not} preserve quasi-coherence in general. To 
remedy this, one uses the \textit{coherator functor} $Q: 
\mathcal{O}_X\text{-Mod} \rightarrow \Q(X)$, which is the right adjoint 
to the inclusion functor $i: \Q(X) \hookrightarrow 
\mathcal{O}_X\text{-Mod}$ (see \cite{PIS}). The internal hom in $\Q(X)$ 
is then defined as
$$\hm(\mathcal{F}, \mathcal{G}) := Q\!\left(\hom(\mathcal{F}, 
\mathcal{G})\right),$$
for any $\mathcal{F}, \mathcal{G} \in \Q(X)$. Together with the tensor 
product $\otimes_{\O}$, this endows $\Q(X)$ with the structure of a 
closed symmetric monoidal Grothendieck category (see \cite{PIS}). In 
particular, for any object $\mathcal{G} \in \Q(X)$, the pair 
$(-\otimes_{\O}\mathcal{G},\ \hm(\mathcal{G},-))$ forms an adjoint pair. 
We fix once and for all an injective cogenerator $I$ of $\Q(X)$, and 
we set $(-)^{+}:= \hm(-,I)$.

\medskip
This monoidal structure allows us to define purity in $\Q(X)$ via the 
tensor product $\otimes_{\mathcal{O}_X}$. In general, there are four 
notions of purity: categorical purity, tensor purity, local purity and 
stalk-wise purity (or geometric purity). The following results, 
focusing on geometric purity in $\Q(X)$, are collected from \cite{PIS} 
and \cite{PMONOIDAL}.
\begin{itemize}
\item A short exact sequence $\mathcal{E}$ in $\Q(X)$ is called 
\textit{pure} if for any object $\mathcal{M}$ in $\Q(X)$, 
$\mathcal{M}\otimes\mathcal{E}$ remains exact.
\item A monomorphism $\mathcal{N}\longrightarrow \mathcal{K}$ in $\Q(X)$ 
is called \textit{pure} if its corresponding short exact sequence is pure.
\item An object $\mathcal{A}$ in $\Q(X)$ is called \textit{absolutely 
pure} if any monomorphism $\mathcal{A}\longrightarrow \mathcal{K}$ in 
$\Q(X)$ is pure.
\item An object $\mathcal{F}$ in $\Q(X)$ is said to be \textit{flat} if 
the functor $\mathcal{F}\otimes- : \Q(X)\longrightarrow\Q(X)$ preserves 
monomorphisms.
\end{itemize}

\par The following conditions hold by \cite{PMONOIDAL}:
\begin{enumerate}
\item For every short exact sequence $\mathcal{E}$ in $\Q(X)$, 
$\mathcal{E}$ is pure if and only if $\mathcal{E}^{+}$ splits.
\item If $X$ is semi-separated, the following conditions are equivalent for any object $\mathcal{F}$ 
in $\Q(X)$:
\begin{itemize}
\item[(i)] $\mathcal{F}$ is flat.
\item[(ii)] $\mathcal{F}^{+}$ is injective in $\Q(X)$.
\item[(iii)] Any short exact sequence ending in $\mathcal{F}$ is pure.
\end{itemize}
\item An object $\mathcal{A}$ in $\Q(X)$ is absolutely pure if and only 
if it is a pure subobject of an injective object.
\end{enumerate}

\par According to \cite[Proposition 4.3 and Corollary 4.6]{PIS}, for 
any $\mathcal{F}$ in $\Q(X)$, $\mathcal{F}^{+}$ is pure injective.
 
\par Recall that, given a class of objects $\mathcal{L}$ in an abelian category $\mathcal{C}$,  we associate its right first Ext-orthogonal class, $$\mathcal{L}^\perp= \{X \in \mathcal{C} \mid \Ext^1(L,X) = 0, L\in \mathcal{L}\},$$ and its left first Ext-orthogonal class, $${}^\perp\mathcal{L}= \{X \in \mathcal{C} \mid \Ext^1(X,L) = 0, L\in \mathcal{L}\}.$$ In particular, if $\mathcal{L}= \{M\}$ then we simply write ${}^\perp\mathcal{L}={}^\perp M$ and  $\mathcal{L}^\perp=M^\perp$.

\par Throughout, we denote by $\Flat(X)$, $\Abs(X)$ and 
$\inj(X)$ the classes of flat, absolutely pure and injective 
objects of $\Q(X)$, respectively.

\section{Absolutely pure domains in $\Q(X)$} \label{sec 2}

\par In this section, we introduce and study absolutely pure domains in $\Q(X)$.

\begin{Def}\label{def 1}
Given two objects $\mathcal{M}$ and $\mathcal{N}$ in $\Q(X)$, $\mathcal{N}$ is said to be \textit{absolutely $\mathcal{M}$-pure} if $\mathcal{M}\otimes \mathcal{N}\rightarrow \mathcal{M}\otimes \mathcal{K}$ is a monomorphism for every monomorphism $\mathcal{N} \to \mathcal{K}$.

The \textit{absolutely pure domain} of $\mathcal{M}$ is defined as the class
$$\ApQ(\mathcal{M}):=\{\mathcal{N} \in \Q(X): \mathcal{N} \ \text{is absolutely}\ \mathcal{M}\text{-pure}\}.$$
\end{Def}
\begin{prop}\label{prop:ap-in}
Let $\mathcal{M}\in\mathfrak{Qcoh}(X)$. Then
$\mathcal{A}p_{\mathfrak{Qcoh}(X)}(\mathcal{M})=\underline{\mathfrak{In}}^{-1}_{\mathfrak{Qcoh}(X)}(\mathcal{M}^{+})$.
\end{prop}

\begin{proof}
Let $g\colon\mathcal{N}\to\mathcal{K}$ be a monomorphism. By the adjunction between
$-\otimes\mathcal{M}$ and $\mathcal{H}om^{qc}_{\mathcal{O}_X}(\mathcal{M},-)$, the map
$\mathrm{Hom}(\mathcal{K},\mathcal{M}^{+})\to\mathrm{Hom}(\mathcal{N},\mathcal{M}^{+})$ induced by $g$
is identified with $\mathrm{Hom}(\mathcal{M}\otimes\mathcal{K},I)\to\mathrm{Hom}(\mathcal{M}\otimes\mathcal{N},I)$.
Since $I$ is an injective cogenerator, this map is surjective if and only if
$\mathcal{M}\otimes\mathcal{N}\to\mathcal{M}\otimes\mathcal{K}$ is a monomorphism.
Hence $\mathcal{M}^{+}$ is $\mathcal{N}$-subinjective if and only if
$\mathcal{N}\in\mathcal{A}p_{\mathfrak{Qcoh}(X)}(\mathcal{M})$.
\end{proof}

\par\label{proprieties} The absolutely pure domain $\ApQ(\mathcal{M})$ satisfies the following 
properties, which generalize \cite[Propositions 2.2--2.5]{Dur}:
\begin{enumerate}
\item $\mathcal{M}$ is flat if and only if $\ApQ(\mathcal{M}) = \Q(X)$.
\item $\mathcal{N} \in \ApQ(\mathcal{M})$ if and only if $\mathcal{M} 
\otimes \mathcal{N} \to \mathcal{M} \otimes E(\mathcal{N})$ is a 
monomorphism, where $E(\mathcal{N})$ is the injective envelope of 
$\mathcal{N}$.

\item $\Abs(X) \subseteq \ApQ(\mathcal{M})$ for any $\mathcal{M} \in 
\Q(X)$, and the intersection of the absolutely pure domains of all 
objects of $\Q(X)$ coincides with the class of absolutely pure objects.

\item For a family 
$(\mathcal{M}_i)_{i \in I}$ of objects in $\Q(X)$,
$$\ApQ\!\left(\bigoplus_{i \in I} \mathcal{M}_i\right) = \bigcap_{i \in I} \ApQ(\mathcal{M}_i),$$

\item If $\mathcal{M}'$ is a pure subobject of $\mathcal{M}$, then
\[
\ApQ(\mathcal{M}) = \ApQ(\mathcal{M}') \cap \ApQ(\mathcal{M}/\mathcal{M}').
\]
\end{enumerate}

The following lemma gives a further characterization of $\ApQ(\mathcal{M})$
in terms of factorization through injective objects, which we use below.

\begin{lem}\label{lem-factor}
Let $\mathcal{M},\mathcal{N}\in\Q(X)$. Then $\mathcal{N}\in\ApQ(\mathcal{M})$ 
if and only if every morphism $\mathcal{N}\to\mathcal{M}^{+}$ factors 
through an injective object of $\Q(X)$.
\end{lem}

\begin{proof}
Let $\iota:\mathcal{N}\rightarrow E(\mathcal{N})$ be the injective envelope
of $\mathcal{N}$. By \ref{proprieties}(2),
\[
\mathcal{N}\in\ApQ(\mathcal{M})
\]
if and only if
\[
\mathcal{M}\otimes\mathcal{N}\longrightarrow
\mathcal{M}\otimes E(\mathcal{N})
\]
is a monomorphism. Since $I$ is an injective cogenerator of $\Q(X)$, this
is equivalent to the surjectivity of
\[
\Hom(\mathcal{M}\otimes E(\mathcal{N}),I)
\longrightarrow
\Hom(\mathcal{M}\otimes\mathcal{N},I).
\]
Using the adjunction
\[
(-\otimes\mathcal{M},\ \hm(\mathcal{M},-)),
\]
this is equivalent to the surjectivity of the restriction map
\[
\Hom(E(\mathcal{N}),\mathcal{M}^{+})
\longrightarrow
\Hom(\mathcal{N},\mathcal{M}^{+}),
\qquad
g\longmapsto g\circ\iota.
\]

Thus, $\mathcal{N}\in\ApQ(\mathcal{M})$ if and only if every morphism
$f:\mathcal{N}\rightarrow\mathcal{M}^{+}$ extends to
$E(\mathcal{N})$. Since $E(\mathcal{N})$ is injective, this is
equivalent to saying that every morphism
$f:\mathcal{N}\rightarrow\mathcal{M}^{+}$ factors through an injective
object. Indeed, one implication is immediate. Conversely, if
$f=h\circ p$, where $p:\mathcal{N}\rightarrow\mathcal{J}$ and
$\mathcal{J}$ is injective, then $p$ extends to a morphism
$\widetilde p:E(\mathcal{N})\rightarrow\mathcal{J}$, and
$h\circ\widetilde p$ is an extension of $f$.
\end{proof}

\par The following result provides a new characterization of the 
absolutely pure domain of some quasi-coherent $\O$-modules via Ext-orthogonal classes. 
\begin{thm}\label{M-F-M}
Let $0\to \mathcal{M} \to \mathcal{F}\to \mathcal{M}'\to 0$ be a short exact 
sequence in $\Q(X)$ with $\mathcal{F}$ flat. Then, 
$^{\perp}{\mathcal{M}'^{+}}\subseteq \ApQ(\mathcal{M})$. If moreover 
$\inj(X)\subseteq {^{\perp}\mathcal{M}'^{+}}$, then 
${^{\perp}\mathcal{M}'^{+}}=\ApQ(\mathcal{M})$.
\end{thm}

\begin{proof}
Let $\mathcal{N}\in{^{\perp}{\mathcal{M}'^{+}}}$. We apply the functor 
$\Hom(\mathcal{N},-)$  to the short exact sequence 
$0\to \mathcal{M}'^{+} \to \mathcal{F}^{+}\to \mathcal{M}^{+}\to 0$, 
resulting in the following long exact sequence:
$$\xymatrix{\ar[r]& \Hom(\mathcal{N},\mathcal{F}^{+}) \ar[r]& 
\Hom(\mathcal{N},\mathcal{M}^{+}) \ar[r]& 
\Ext^1(\mathcal{N},\mathcal{M}'^{+})\ar[r]& 
\Ext^1(\mathcal{N},\mathcal{F}^{+}) \ar[r]&}$$
Since $\mathcal{N}\in{^{\perp}{\mathcal{M}'^{+}}}$, 
$\Ext^1(\mathcal{N},\mathcal{M}'^{+})=0$, which implies that 
$\Hom(\mathcal{N},\mathcal{F}^{+}) \to \Hom(\mathcal{N},\mathcal{M}^{+})$ 
is surjective. Thus, any morphism $\mathcal{N}\to \mathcal{M}^{+}$ factors 
through the injective object $\mathcal{F}^{+}$, so 
$\mathcal{N}\in \ApQ(\mathcal{M})$.

\par Additionally, suppose that 
$\inj(X))\subseteq {^{\perp}{\mathcal{M}'^{+}}}$. Let 
$\mathcal{N}\in \ApQ(\mathcal{M})$, and let $\mathcal{N}\to \mathcal{E}$ be a 
monomorphism with $\mathcal{E}$ injective. We apply the functors $\Hom(\mathcal{E},-)$ and 
$\Hom(\mathcal{N},-)$ to the morphism $\mathcal{F}^{+}\to \mathcal{M}^{+}$, 
resulting in the following commutative diagram:
$$
\xymatrix{
\Hom(\mathcal{E},\mathcal{F}^{+})\ar[r]\ar[d]& 
\Hom(\mathcal{E},\mathcal{M}^{+})\ar[d]\ar[r] & 
\Ext^1(\mathcal{E},\mathcal{M}'^{+}) \ar[r]&\Ext^1(\mathcal{E},\mathcal{F}^{+})\\
\Hom(\mathcal{N},\mathcal{F}^{+}) \ar[r]& 
\Hom(\mathcal{N},\mathcal{M}^{+})\ar[r] & 
\Ext^1(\mathcal{N},\mathcal{M}'^{+}) \ar[r]& 
\Ext^1(\mathcal{N},\mathcal{F}^{+})\\
}
$$
Since $\mathcal{F}^{+}$ is injective, $\Ext^1(\mathcal{N},\mathcal{F}^{+})=0$. 
Since $\mathcal{N}\in \ApQ(\mathcal{M})$, the morphism 
$\Hom(\mathcal{E},\mathcal{M}^{+}) \to \Hom(\mathcal{N},\mathcal{M}^{+})$ is 
surjective, and $\Hom(\mathcal{E},\mathcal{F}^{+})\to\Hom(\mathcal{E},\mathcal{M}^{+})$ is 
surjective (as $\Ext^1(\mathcal{E},\mathcal{M}'^{+})=0$ by assumption). This implies 
that $\Hom(\mathcal{N},\mathcal{F}^{+})\to \Hom(\mathcal{N},\mathcal{M}^{+})$ 
is surjective, hence $\Ext^1(\mathcal{N},\mathcal{M}'^{+})=0$, so 
$\mathcal{N}\in{^{\perp}\mathcal{M}'^{+}}$.
\end{proof}
 Proposition \ref{A} and Proposition \ref{B} demonstrate that in the case where $X=\spec(A)$ and $\widetilde{N}, \widetilde{M} \in \Q(X)$, one of the two assertions:
\begin{itemize}
	\item $\widetilde{N}|{{D_A}(f_i)}\in\Apf(\widetilde{M}|{{D_A}(f_i)})$, for $i=1,...,n$, where $f_1,...,f_n\in A$ and $X=\cup_{i=1}^{i=n}\D_A(f_i)$.
	\item $\widetilde{N}_{x}\in\Apxx(\widetilde{M}_x)$ for every $x\in X$
\end{itemize} implies that $\widetilde{N}\in \Apq(\widetilde{M})$.

\par The following two lemmas are usefull to prove Proposition \ref{A} and Proposition \ref{B}.

\begin{lem}\label{lem-module-1}
Let $M$ and $N$ be $A$-modules, and $f_1,...,f_n$ be elements generating the unit ideal. If $N_{f_i}$ is absolutely $M_{f_i}$-pure over $A_{f_i}$ for every $i = 1, . . . , n$, then $N$ is absolutely $M$-pure over $A$.
\end{lem}

\begin{proof}
Let $N \to S$ be a monomorphism, and let $K = \Ker(M\otimes N\to M\otimes S)$. For every $i \in \{1,...,n\}$, $K_{f_i}\cong \Ker(M_{f_i}\otimes N_{f_i}\to M_{f_i}\otimes S_{f_i})=0$ since $N_{f_i}$ is absolutely $M_{f_i}$-pure over $A_{f_i}$ by assumption. Thus, for every $x\in K$ and any $i\in\{1,...,n\}$, there exists some $s_i\geq 0$ such that $f_i^{s_i}x=0$.
On the other hand, since $1=h_1f_1+...+h_nf_n$ for some $h_1,...,h_n\in A$, for any $x\in K$, we can write $x = 1.x = 1^{s}x = (h_1f_1+...+h_nf_n)^{s}x=0$, where $s > s_1+...+s_n-1$. Thus, $K=0$, and $N$ is absolutely $M$-pure.
\end{proof}

\begin{lem}\label{lem-module-2}
Let $M$ and $N$ be $A$-modules. If $N_x$ is absolutely $M_x$-pure over $A_x$ for every prime ideal $x$, then $N$ is absolutely $M$-pure over $A$.
\end{lem}
\begin{proof}
Let $N \to K$ be a monomorphism. For any $x\in X$, $M_x\otimes N_x\rightarrow M_x\otimes K_x$ is monic by hypothesis. The tensor product functor commutes with the stalk functor and reflects monomorphisms. Therefore, $M\otimes N\rightarrow M\otimes K$ is injective, so $N$ is absolutely $M$-pure over $A$.
\end{proof}

\begin{prop}\label{A}
Let $X=\spec(A)$, $f_1,...,f_n\in A$ such that $X=\cup_{i=1}^{i=n}\D_A(f_i)$, and $\widetilde{N},\widetilde{M}\in \Q(X)$. If $\widetilde{N}|_{{D_A}(f_i)}\in\ApQ(\widetilde{M}|_{{D_A}(f_i)})$ for any $i\in\{1,...,n\}$, then $\widetilde{N}\in \ApQ(\widetilde{M})$.
\end{prop}
\begin{proof}
The category of quasi-coherent $\O$-modules is equivalent to the category of modules over $A$. Therefore, by Lemma \ref{lem-module-1}, the result follows.
\end{proof}

\begin{prop}\label{B}
Let $X=\spec(A)$, and $\widetilde{N}, \widetilde{M} \in \Q(X)$. If $\widetilde{N}_{x}\in\Apxx(\widetilde{M}_x)$ for every $x\in X$, then $\widetilde{N}\in \ApQ(\widetilde{M})$.
\end{prop}
\begin{proof}
The result follows from Lemma \ref{lem-module-2}.
\end{proof}

\par In general, the converse statements of Propositions \ref{A} and \ref{B} are not valid (see Remark \ref{AA}). Nevertheless, when $A$ is coherent, these converse statements do hold, as we will see in Proposition \ref{help coh}. Consequently, it is reasonable to introduce and explore the notion of locally absolutely pure domains for quasi-coherent sheaves on any scheme, particularly on locally coherent schemes.

\begin{Def}
Let $\mathcal{M}$ and $\mathcal{N}$ be quasi-coherent $\O$-modules. We say that $\mathcal{N}$ is \textit{locally absolutely $\mathcal{M}$-pure} if $\mathcal{N}(U)\in\Apou(\mathcal{M}(U))$ for every affine open subset $U\subseteq X$.

The \textit{locally absolutely pure domain} of $\mathcal{M}$ is defined as the class
$$\LAp(\mathcal{M}):=\{\mathcal{N} \in \Q(X): \mathcal{N} \ \text{is locally absolutely}\ \mathcal{M}\text{-pure}\}.$$
\end{Def}

\par The following proposition  states that, to establish whether a quasi-coherent $\O$-module  $\mathcal{N}$ is locally absolutely $\mathcal{M}$-pure over a locally coherent scheme, it suffices to examine a cover consisting of affine open subsets of $X$. Furthermore, this local property holds stalkwise.

\begin{prop}\label{locally}
Let $X$ be a locally coherent scheme and $\mathcal{M},\mathcal{N}\in\Q(X)$. Then, the following conditions are equivalent:
\begin{itemize}
\item[(i)] $\mathcal{N} \in \LAp(\mathcal{M})$.
\item[(ii)] For any cover $\mathcal{U}= \lbrace U_i\rbrace_{i\in I}$ of affine open subsets of $X$, $\mathcal{N}(U_i)\in \Apoui(\mathcal{M}(U_i))$ for any $i\in I$.
\item[(iii)] There exists a cover $\mathcal{U}= \lbrace U_i\rbrace_{i\in I}$ of affine open subsets of $X$ such that, for any $i\in I$, $\mathcal{N}(U_i)\in \Apoui(\mathcal{M}(U_i))$.
\item[(iv)] $\mathcal{N}_{x}\in \Apxx(\mathcal{M}_x)$ for all $x\in X$.
\end{itemize}
\end{prop}
\begin{proof}
The implications $(i)\Rightarrow (ii)\Rightarrow (iii)$ are straightforward.
\par $(iii)\Rightarrow (iv)$: Let $U_i\in \mathcal{U}$ and $\mathcal{N}(U_i)\rightarrow E$ be an injective envelope in the category of $\mathcal{O}(U_i)$-modules. By assumption, $\mathcal{M}(U_i)\otimes\mathcal{N}(U_i)\rightarrow \mathcal{M}(U_i)\otimes E$, which implies that $\mathcal{M}|_{U_i}\otimes\mathcal{N}|_{U_i}\rightarrow \mathcal{M}|_{U_i}\otimes \widetilde{E}$ is monic in $\Q(U_i)$. Thus, for any $x\in U_i$, $\mathcal{M}_x\otimes\mathcal{N}_x\rightarrow \mathcal{M}_x\otimes E_x$ is also monic. Since $E_x$ is absolutely pure (as shown in \cite[Theorem 3.21]{Ab} for $\mathcal{O}_{X,x}$ being coherent), then $\mathcal{N}_{x}\in \Apxx(\mathcal{M}_x)$ for any $x\in U_i$. Since $\mathcal{U}$ is a cover of $X$, we have $\mathcal{N}_{x}\in \Apxx(\mathcal{M}_x)$ for any $x\in X$.
\par $(iv)\Rightarrow (i)$: Let $U$ be an affine open subset of $X$, and consider a monomorphism $\mathcal{N}(U)\rightarrow K$ in the category of $\mathcal{O}(U)$-modules. Thus, $\mathcal{N}|_{U}\rightarrow \widetilde{K}$ is monic in $\Q(U)$. For any $x\in U$, $\mathcal{M}_x \otimes \mathcal{N}_x\rightarrow \mathcal{M}_x\otimes K_x$ is a monomorphism since $\mathcal{N}_x\rightarrow K_x$ is monic and $\mathcal{N}_{x}\in \Apxx(\mathcal{M}_x)$. Thus, $\mathcal{M}|_{U} \otimes \mathcal{N}|_{U}\rightarrow \mathcal{M}|_{U}\otimes\widetilde{K}$ is a monomorphism, which implies that $\mathcal{M}(U) \otimes \mathcal{N}(U)\rightarrow \mathcal{M}(U)\otimes K$ is monic. Therefore, $\mathcal{N}$ is locally absolutely $\mathcal{M}$-pure.
\end{proof}

\begin{thm}\label{help coh}
Let $X$ be a locally coherent scheme. Then, for every quasi-coherent $\O$-module $\mathcal{M}$, $\LAp(\mathcal{M}) \subseteq \ApQ(\mathcal{M})$. If moreover $X$ is semi-separated and quasi-compact, then $\LAp(\mathcal{M}) = \ApQ(\mathcal{M})$.
\end{thm}
\begin{proof}
For any $\mathcal{N} \in \LAp(\mathcal{M})$, Propositions \ref{A}, \ref{B} and \ref{locally} imply that $\mathcal{N}\in\ApQ(\mathcal{M})$.
\par Assuming that $X$ is semi-separated and quasi-compact, and $\mathcal{N}\in \ApQ(\mathcal{M})$, we aim to prove that $\mathcal{N}$ is locally absolutely $\mathcal{M}$-pure. By Proposition \ref{locally}, it suffices to show that $\mathcal{N}_x\in \Apxx(\mathcal{M}_x)$ for every $x\in X$. Consider $\mathcal{N}\rightarrow \mathcal{E}$ as an injective envelope in $\Q(X)$. Then, $\mathcal{M}\otimes\mathcal{N}\rightarrow\mathcal{M}\otimes\mathcal{E}$ is monic, which implies that $\mathcal{M}_x\otimes\mathcal{N}_x\rightarrow\mathcal{M}_x\otimes\mathcal{E}_x$ for every $x\in X$. Since $X$ is locally coherent, semi-separated and quasi-compact, by \cite[Proposition 3.3]{absolute}, $\mathcal{E}_x$ is absolutely pure for every $x\in X$. Thus, $\mathcal{N}_x\in \Apxx(\mathcal{M}_x)$ for every $x\in X$ by Proposition \ref{def 1}.
\end{proof}

\begin{thm}\label{coh}
Let $X$ be a semi-separated quasi-compact scheme. The following are equivalent:
\begin{itemize}
\item[(i)] $X$ is locally coherent.
\item[(ii)] For every quasi-coherent $\O$-module $\mathcal{M}$, $\ApQ(\mathcal{M}) \subseteq \LAp(\mathcal{M})$.
\end{itemize}
\end{thm}
\begin{proof}
$(i)\Rightarrow (ii)$: Proposition \ref{help coh}.
\par $(ii)\Rightarrow (i)$: We show that any absolutely pure quasi-coherent sheaf is locally absolutely pure. Let $U$ be an affine open subset of $X$ and $i: U \rightarrow X$ be the
inclusion. Since $X$ is semi-separated, the direct image functor $i_{U,*}$ preserves quasi-coherence. Let $\mathcal{N}$ be an absolutely pure quasi-coherent $\O$-module, and $U$ an affine open subset of $X$. For any $\O(U)$-module $M$, $\mathcal{N}$ is absolutely $i_{U,*}\widetilde{M}$-pure. Since $i_{U,*}\widetilde{M}(U)=M$, by hypothesis, $\mathcal{N}(U)$ is absolutely $M$-pure, meaning that $\mathcal{N}(U)$ is absolutely pure in $\XU$. As a result, $\mathcal{N}$ is locally absolutely pure.
\end{proof}

\section{$F$-indigent objects in $\Q(X)$}

\par In this section, we introduce and study $f$-indigent objects in $\Q(X)$.

Durğun in \cite{Dur} explored the concept of flat indigent modules as the opposite concept of flat modules. Here, we extend this notion to $\Q(X)$.

\begin{Def}\label{def f-indigent}
An object $\mathcal{M}$ in $\Q(X)$ is said to be \textit{flat indigent} (or an \textit{$f$-indigent object}) if $\ApQ(\mathcal{M})=\{\mathcal{N}\in\Q(X):\, \mathcal{N}\text{ is absolutely pure}\}$.
\end{Def}
\par Recall that the category $\Q(X)$ on a quasi-compact and quasi-separated scheme is locally finitely presented (see \cite[Corollary 6.9.12]{alexander}).
\begin{lem}\label{presented}
Assume that $X$ is quasi-compact and quasi-separated. Then, an object $\mathcal{N}$ in $\Q(X)$ is absolutely pure if and only if $\mathcal{M}\otimes\mathcal{N}\rightarrow\mathcal{M}\otimes E(\mathcal{N})$ is a monomorphism for each finitely presented object $\mathcal{M}$.
\end{lem}
\begin{proof}
The direct implication is clear.
\par Conversely, let $\mathcal{G}$ be an object in $\Q(X)$, $\mathcal{G}=\varinjlim \mathcal{M}_i$, where $\{\mathcal{M}_i\}_{i\in I}$ is a direct system of finitely presented objects in $\Q(X)$. By assumption, we have the exact sequences $\mathcal{E}_i: 0\rightarrow\mathcal{M}_i\otimes\mathcal{N}\rightarrow\mathcal{M}_i\otimes E(\mathcal{N})$ for all $i\in I$. Since the functor direct limit is exact and commutes with tensor product, we get the exact sequence $\varinjlim \mathcal{E}_i:0\rightarrow\mathcal{G}\otimes\mathcal{N}\rightarrow\mathcal{G}\otimes E(\mathcal{N})$. So, $\mathcal{N}$ is absolutely pure.
\end{proof}
\par The following result generalizes \cite[Proposition 3.1]{Dur}, 
which states that every ring has an $f$-indigent module, to the 
category $\Q(X)$.

\begin{prop}\label{f-indigent}
Assume that $X$ is quasi-compact and quasi-separated. Then, $\Q(X)$ has 
an $f$-indigent object.
\end{prop}

\begin{proof}
Let $\Lambda$ be a complete set of representatives of finitely presented 
objects in $\Q(X)$. Set $\mathcal{G}:=\bigoplus_{\mathcal{T}_i\in 
\Lambda}\mathcal{T}_i$. We claim that $\mathcal{G}$ is an $f$-indigent 
object. Let $\mathcal{N}\in\ApQ(\mathcal{G})$. By the properties of 
$\ApQ$, we have 
$$\ApQ\!\left(\bigoplus_{\mathcal{T}_i\in\Lambda}\mathcal{T}_i\right) = 
\bigcap_{\mathcal{T}_i\in\Lambda}\ApQ(\mathcal{T}_i),$$
so $\mathcal{N}\in\ApQ(\mathcal{T}_i)$ for every finitely presented object 
$\mathcal{T}_i\in\Lambda$. This means that $\mathcal{T}_i\otimes\mathcal{N}
\rightarrow\mathcal{T}_i\otimes E(\mathcal{N})$ is a monomorphism for every 
finitely presented object $\mathcal{T}_i$. Therefore, by Lemma 
\ref{presented}, $\mathcal{N}$ is absolutely pure in $\Q(X)$. Hence, 
$\ApQ(\mathcal{G})\subseteq\Abs(X)$. Since the reverse inclusion 
$\Abs(X)\subseteq\ApQ(\mathcal{G})$ always holds, we conclude that 
$\ApQ(\mathcal{G})=\Abs(X)$, and thus $\mathcal{G}$ is $f$-indigent.
\end{proof}
\par The following propositions are direct consequences of the properties of $\ApQ$ established in Section \ref{proprieties}.

\begin{prop}
If an object in $\Q(X)$ has an $f$-indigent pure subobject, then it is $f$-indigent itself.
\end{prop}

\begin{prop}
Let $\mathcal{P}$ be a subobject of $\mathcal{M}$. If $\mathcal{M}/\mathcal{P}$ is flat, then $\mathcal{P}$ is $f$-indigent if and only if $\mathcal{M}$ is $f$-indigent.
\end{prop}

\begin{proof}
Since $\mathcal{M}/\mathcal{P}$ is flat, $\mathcal{P}$ is pure in $\mathcal{M}$ and $\ApQ(\mathcal{M}/\mathcal{P})=\Q(X)$. Then $\ApQ(\mathcal{P})=\ApQ(\mathcal{P})\cap \ApQ(\mathcal{M}/\mathcal{P})= \ApQ(\mathcal{M})$, which implies that $\mathcal{P}$ is $f$-indigent if and only if $\mathcal{M}$ is $f$-indigent.
\end{proof}

\begin{prop}\label{prop:abs-inj}
Let $X$ be a scheme.
\begin{enumerate}
\item[(a)] If $\mathrm{Qcoh}(X)$ has a pure injective indigent object $\mathcal{M}$, then
$\mathrm{Abs}(X)=\mathrm{Inj}(X)$ and $\mathcal{M}^{+}$ is a pure injective
$f$-indigent object.
\item[(b)] If $X$ is quasi-compact and quasi-separated and
$\mathrm{Abs}(X)=\mathrm{Inj}(X)$, then $\mathrm{Qcoh}(X)$ has a pure
injective indigent object.
\end{enumerate}
In particular, for $X$ quasi-compact and quasi-separated, the two conditions are equivalent.
\end{prop}

\begin{proof}
(a) Let $\mathcal{M}$ be a pure injective indigent object of $\mathrm{Qcoh}(X)$. The evaluation
morphism $\nu_{\mathcal{M}}\colon\mathcal{M}\to\mathcal{M}^{++}$ is a pure monomorphism, hence
splits since $\mathcal{M}$ is pure injective: $\mathcal{M}^{++}\cong\mathcal{M}\oplus\mathcal{C}$
with $\mathcal{C}=\mathcal{M}^{++}/\mathcal{M}$. Therefore
\[
\mathrm{In}^{-1}_{\mathrm{Qcoh}(X)}(\mathcal{M}^{++})
=\mathrm{In}^{-1}{\mathrm{Qcoh}(X)}(\mathcal{M})\cap\mathrm{In}^{-1}{\mathrm{Qcoh}(X)}(\mathcal{C})
\subseteq\mathrm{In}^{-1}_{\mathrm{Qcoh}(X)}(\mathcal{M})=\mathrm{Inj}(X).
\]
By Proposition~\ref{prop:ap-in} applied to $\mathcal{M}^{+}$ (note $(\mathcal{M}^{+})^{+}=\mathcal{M}^{++}$),
\[
\mathrm{Ap}_{\mathrm{Qcoh}(X)}(\mathcal{M}^{+})
=\mathrm{In}^{-1}_{\mathrm{Qcoh}(X)}(\mathcal{M}^{++})\subseteq\mathrm{Inj}(X).
\]
Since $\mathrm{Inj}(X)\subseteq\mathrm{Abs}(\mathrm{Qcoh}(X))\subseteq\mathrm{Ap}_{\mathrm{Qcoh}(X)}(\mathcal{M}^{+})$
always holds, all these classes coincide. Hence $\mathrm{Abs}(\mathrm{Qcoh}(X))=\mathrm{Inj}(X)$ and
$\mathrm{Ap}_{\mathrm{Qcoh}(X)}(\mathcal{M}^{+})=\mathrm{Abs}(\mathrm{Qcoh}(X))$, i.e.\
$\mathcal{M}^{+}$ is $f$-indigent; it is pure injective, being a dual object.

(b) Since $X$ is quasi-compact and quasi-separated, $\mathrm{Qcoh}(X)$ is locally finitely
presented. Let $\Lambda$ be a complete set of representatives of finitely presented objects.
By Proposition~\ref{f-indigent}, $\mathcal{G}:=\bigoplus_{\mathcal{T}_i\in\Lambda}\mathcal{T}_i$
is $f$-indigent. By Proposition\ref{prop:ap-in},
\[
\mathrm{In}^{-1}_{\mathrm{Qcoh}(X)}(\mathcal{G}^{+})
=\mathrm{Ap}_{\mathrm{Qcoh}(X)}(\mathcal{G})
=\mathrm{Abs}(\mathrm{Qcoh}(X))=\mathrm{Inj}(X).
\]
Thus $\mathcal{G}^{+}$, which is pure injective, is indigent.
\end{proof}

\begin{prop}\label{f-sum}
If $\mathcal{M}$ is an $f$-indigent object of $\Q(X)$, then $\mathcal{M}\oplus \mathcal{G}$ is $f$-indigent for any object $\mathcal{G}\in \Q(X)$.
\end{prop}

\begin{proof}
Since $\ApQ(\mathcal{M}\oplus \mathcal{G})=\ApQ(\mathcal{M})\cap\ApQ(\mathcal{G})$ and $\ApQ(\mathcal{M})$ consists only of absolutely pure objects, then $\ApQ(\mathcal{M}\oplus \mathcal{G})=\ApQ(\mathcal{M})$. Hence, $\mathcal{M}\oplus \mathcal{G}$ is $f$-indigent.
\end{proof}

\begin{rem}\label{AA}
The converse statements of Propositions \ref{A} and \ref{B} are 
not valid in general, showing that $f$-indigency is not a local 
property in $\Q(X)$. Indeed, let $X = \spec(A)$ and let $M$ be 
an $f$-indigent $A$-module. If the converse of Proposition \ref{A} 
were true, then $f$-indigency would be a local property, meaning 
that for every $\widetilde{N} \in \Apq(\widetilde{M})$, $N$ would be absolutely 
pure over $A$, and $N_{f_i}$ would be absolutely pure over $A_{f_i}$ 
for each $i = 1,...,n$. However, this does not hold without the coherence assumption on $A$.
\end{rem}

\par This failure of locality motivates the introduction of the 
following notion of locally $f$-indigent quasi-coherent modules.

\begin{Def}
Let $\mathcal{M}$ be an object in $\Q(X)$. $\mathcal{M}$ is said to be \textit{locally $f$-indigent} if $\mathcal{M}(U)$ is  $f$-indigent $\O(U)$-module  for every affine open $U\subseteq X$.
\end{Def}

\par The next proposition states that to check whether a quasi-coherent 
$\mathcal{O}_X$-module on a locally coherent scheme is locally 
$f$-indigent, it suffices to look at some cover by affine open subsets 
of $X$. We first recall the skyscraper functor that we need to 
demonstrate this proposition:

If $x\in X$ and $M$ is an 
$\mathcal{O}_{X,x}$-module, then the skyscraper sheaf $i_{x,*}(M)$ is 
an $\mathcal{O}_X$-module defined as:

\[
i_{x,*}(M)(U)=
\begin{cases}
    M & \text{if } x\in U \\
    0 & \text{otherwise}
\end{cases}
\qquad \forall\, U\subseteq X
\]

The stalks of $i_{x,*}(M)$ are $\lbrace 0\rbrace$ for all points except 
at $x$, where $(i_{x,*}(M))_x= M$.\\
This construction yields a functor $i_{x,*}: \Xx\rightarrow 
\mathcal{O}_X\text{-Mod}$, which is fully faithful, left exact, 
preserves injectives and serves as a right adjoint to taking stalks at 
$x$.

\begin{prop}\label{loc ind}
Let $X$ be a locally coherent scheme, and $\mathcal{M}$ be a quasi-coherent $\O$-module. The following assertions are equivalent:
\begin{itemize}
\item[(i)] $\mathcal{M}(U)$ is $f$-indigent for every affine open $U\subseteq X$.
\item[(ii)] $\mathcal{M}(U_i)$ is $f$-indigent for all $i\in I$ for some cover $\mathcal{U}=\lbrace U_i\rbrace_{i\in I}$ of affine open subsets.
\item[(iii)] $\mathcal{M}_{x}$ is $f$-indigent for all $x\in X$.
\end{itemize}
\end{prop}
\begin{proof}
$(i)\Rightarrow (ii)$ is obvious.\\
$(ii)\Rightarrow (iii)$: For a given $x\in X$, there exists an 
affine open subset $U_i \in \mathcal{U}$ such that $x \in U_i$. 
Let $N\in\Apxx(\mathcal{M}_x)$, and let $N \to E$ be an injective 
envelope in $\mathcal{O}_{X,x}$-Mod. Consider the skyscraper sheaf 
$i_{x,*}:\Xx\to\Ui$. Since $i_{x,*}$ is left exact, the monomorphism 
$N \to E$ in $\mathcal{O}_{X,x}$-Mod induces a monomorphism 
$i_{x,*}N \to i_{x,*}E$ in $\mathcal{O}_{U_i}$-Mod.

 We show that $i_{x,*}N\otimes_{\mathcal{O}_{U_i}} 
\mathcal{M}|_{U_i} \to i_{x,*}E\otimes_{\mathcal{O}_{U_i}} 
\mathcal{M}|_{U_i}$ is monic in $\mathcal{O}_{U_i}$-Mod. Since 
the stalk functor commutes with tensor products and detects 
monomorphisms in $\mathcal{O}_{U_i}$-Mod, it suffices to show 
that for every $y \in U_i$:
$$(i_{x,*}N\otimes\mathcal{M}|_{U_i})_y \to 
(i_{x,*}E\otimes\mathcal{M}|_{U_i})_y$$
is monic, i.e. $(i_{x,*}N)_y\otimes_{\mathcal{O}_{X,y}}\mathcal{M}_y 
\to (i_{x,*}E)_y\otimes_{\mathcal{O}_{X,y}}\mathcal{M}_y$ is monic. 
By definition of the skyscraper functor:
\begin{itemize}
\item If $y = x$: $(i_{x,*}N)_x = N$ and $(i_{x,*}E)_x = E$, 
so the map is $\mathcal{M}_x \otimes N \to \mathcal{M}_x \otimes E$, 
which is monic since $N\in\Apxx(\mathcal{M}_x)$ and $N\to E$ is 
monic in $\mathcal{O}_{X,x}$-Mod.
\item If $y \neq x$: $(i_{x,*}N)_y = 0$ and $(i_{x,*}E)_y = 0$, 
so the map $0 \to 0$ is trivially monic.
\end{itemize}
Therefore, $i_{x,*}N\otimes\mathcal{M}|_{U_i} \to 
i_{x,*}E\otimes\mathcal{M}|_{U_i}$ is monic in 
$\mathcal{O}_{U_i}$-Mod. So,
$$i_{x,*}N(U_i)\otimes_{\mathcal{O}_X(U_i)}\mathcal{M}(U_i) \to 
i_{x,*}E(U_i)\otimes_{\mathcal{O}_X(U_i)}\mathcal{M}(U_i)$$
becomes $N\otimes_{\mathcal{O}_X(U_i)}\mathcal{M}(U_i) \to 
E\otimes_{\mathcal{O}_X(U_i)}\mathcal{M}(U_i)$, which is monic.

\par Since $i_{x,*}$ preserves injective objects, $i_{x,*}E$ is injective in 
$\mathcal{O}_{U_i}$-Mod. By the adjoint pair $(\widetilde{-}, 
\Gamma(U_i,-))$ and the exactness of $\widetilde{-}$, the functor 
$\Gamma(U_i,-)$ preserves injective objects (see 
\cite[Proposition~25]{daniel}). Therefore, 
$E = \Gamma(U_i, i_{x,*}E) = i_{x,*}E(U_i)$ is injective in 
$\mathcal{O}_X(U_i)$-Mod. By the properties of $\Apq$, since 
$N \to E$ is a monomorphism with $E$ injective in 
$\mathcal{O}_X(U_i)$-Mod and 
$\mathcal{M}(U_i)\otimes N \to \mathcal{M}(U_i)\otimes E$ is monic, 
we conclude that $N\in\Apoui(\mathcal{M}(U_i))$. Since 
$\mathcal{M}(U_i)$ is $f$-indigent by hypothesis, $N$ is absolutely 
pure as an $\mathcal{O}_X(U_i)$-module. Since $\mathcal{O}_X(U_i)$ 
is coherent, by \cite[Theorem~3.21]{Ab}, $N$ is absolutely pure as 
an $\mathcal{O}_{X,x}$-module, hence $\mathcal{M}_x$ is $f$-indigent.

\par $(iii)\Rightarrow (i)$: Let $U$ be an affine open of $X$, $N\in\Apou(\mathcal{M}(U))$, and $N\rightarrow E$ be an injective envelope in $\XU$. Then, $\mathcal{M}(U)\otimes N\rightarrow \mathcal{M}(U)\otimes E$ is monic, which implies that $\mathcal{M}|_{U}\otimes \widetilde{N}\rightarrow \mathcal{M}|_{U}\otimes \widetilde{E}$ is also monic. So, for every $x\in U$, $\mathcal{M}_x\otimes N_x\rightarrow \mathcal{M}_x\otimes E_x$ is monic. Since $\O(U)$ is coherent, $E_x$ is absolutely pure for every $x\in U$. So, by Proposition \ref{def 1}, $N_x\in \Apxx(\mathcal{M}_x)$ for every $x\in U$. By assumption, $N_x$ is absolutely pure for every $x\in U$, hence $N$ is absolutely pure.
\end{proof}

\begin{cor}\label{help neo}
Let $X$ be a semi-separated quasi-compact locally coherent scheme and $\mathcal{M}$ be a quasi-coherent $\O$-module. Then, $\mathcal{M}$ is $f$-indigent in $\Q(X)$ if $\mathcal{M}$ is locally $f$-indigent.
\end{cor}
\begin{proof}
Let $\mathcal{N}$ be an object of $\ApQ(\mathcal{M})$. Then, by Proposition \ref{help coh}, $\mathcal{N}_x\in\Apxx(\mathcal{M}_x)$. By hypothesis and Proposition \ref{loc ind}, $\mathcal{M}_x$ is $f$-indigent for every $x\in X$, so $\mathcal{N}_x$ is absolutely pure for every $x\in X$. Thus, $\mathcal{N}$ is absolutely pure.
\end{proof}

\begin{rem}\label{inj,abs}
On a locally Noetherian scheme $X$, all classes of locally absolutely pure, absolutely pure, locally injective and injective quasi-coherent sheaves are equal, by \cite[II, proposition 7.17, theorem 7.18]{hartshorne}.
\end{rem}

\begin{thm}\label{noetherian}
Let $X$ be a quasi-compact locally coherent scheme. Then, the following assertions are equivalent:
\begin{itemize}
\item[(i)] $X$ is a Noetherian scheme.
\item[(ii)] There exists a pure injective indigent quasi-coherent $\mathcal{O}_U$-module for every affine open $U$ of $X$.
\item[(iii)] There exists a pure injective indigent quasi-coherent $\mathcal{O}_{U_i}$-module for all $i\in I$ for some cover $\lbrace U_i\rbrace_{i\in I}$ of affine open subsets of $X$.
\end{itemize}
\end{thm}
\begin{proof}
$(i)\Rightarrow (ii)$: Let $U$ be an affine open subset of $X$. Since $X$ is Noetherian, $U$ is also Noetherian (see \cite[Corollary 3.22]{gortz}). The result follows from Proposition \ref{prop:abs-inj} and Remark \ref{inj,abs}.
\par $(ii)\Rightarrow (iii)$: Obvious.
\par $(iii)\Rightarrow (i)$: Let $M$ be an absolutely pure $\O(U_i)$-module. Since $\O(U_i)$ is coherent, $M_x$ is absolutely pure for every $x\in U_i$. Let $\widetilde{M}$ be the associated quasi-coherent $\mathcal{O}_{U_i}$-module. Since $\widetilde{M}_x\cong M_x$, by \cite[Lemma 5.4.]{PIS}, $\widetilde{M}$ is absolutely pure. By hypothesis, $\Q(U_i)$ has a pure injective indigent object, so the class of absolutely pure objects equals the class of injective objects in $\Q(U_i)$. Thus, $\widetilde{M}$ is injective, so $M$ is injective, and $X$ is locally Noetherian. Since $X$ is quasi-compact, $X$ is Noetherian.
\end{proof}

\end{document}